\documentclass{amsart}

\usepackage{amsmath,amsfonts,amsthm,amssymb,verbatim}

\newtheorem{proposition}[equation]{Proposition}
\newtheorem{corollary}[equation]{Corollary}
\newtheorem{theorem}[equation]{Theorem}
\newtheorem{conjecture}[equation]{Conjecture}

\numberwithin{equation}{section}

\begin{document}

\title[$p$-adic estimates for multiplicative character sums]{$p$-adic estimates for \\ multiplicative character sums}
\author{Alan Adolphson}
\address{Department of Mathematics\\ Oklahoma State University\\ Stillwater, OK 74078}
\email{adolphs@math.okstate.edu}
\author{Steven Sperber}
\address{School of Mathematics\\ University of Minnesota\\ Minneapolis, MN 55455}
\email{sperber@math.umn.edu}
\date{\today}
\begin{abstract}
This article is an expanded version of the talk given by the first author at the conference ``Exponential sums over finite fields and applications'' (ETH, Z\"urich, November, 2010).  
We state some conjectures on archimedian and $p$-adic estimates for multiplicative character sums over smooth projective varieties.  We also review some of the results of J. Dollarhide\cite{D}, which formed the basis for these conjectures.  Applying his results, we prove one of the conjectures when the smooth projective variety is ${\mathbb P}^n$ itself.
\end{abstract}
\maketitle

\section{Introduction}

Let ${\mathbb F}_q$ be a finite field of characteristic $p$, let $f_1,\dots,f_r\in {\mathbb F}_q[x_0,\dots,x_n]$ be nonconstant homogeneous polynomials of degrees $d_1,\dots,d_r$, respectively, let 
\[ \chi_1,\dots,\chi_r:{\mathbb F}_q^\times\to{\mathbb Q}(\zeta_{q-1})^\times \]
be multiplicative characters (extended to ${\mathbb F}_q$ by setting $\chi_i(0) = 0$ for all $i$), and let $X\subseteq{\mathbb P}^n$ be a projective variety.  We always impose the hypothesis that
\begin{equation}
\text{the product $\chi_1^{d_1}\cdots\chi_r^{d_r}$ is the trivial character.}
\end{equation}
This guarantees that for any point $x\in{\mathbb P}^n({\mathbb F}_q)$, the expression
\[ \chi_1(f_1(x))\cdots \chi_r(f_r(x)) \]
is well-defined, i.e., it is independent of the choice of homogeneous coordinates for the point $x$.  For any $m\geq 1$, we can then form the multiplicative character sums
\begin{equation}
S_m = \sum_{x\in X({\mathbb F}_{q^m})} \prod_{i=1}^r \chi_i\circ{\rm Norm}_{{\mathbb F}_{q^m}/{\mathbb F}_q}(f_i(x))\in{\mathbb Q}(\zeta_{q-1}) \end{equation}
and their associated $L$-function
\begin{equation}
L(t) = \exp\bigg(\sum_{m=1}^\infty S_m\frac{t^m}{m}\bigg)\in{\mathbb Q}(\zeta_{q-1})[[t]]. 
\end{equation}
It is known that this series is a rational function, i.e., $L(t)\in{\mathbb Q}(\zeta_{q-1})(t)$, and one wants to determine the archimedian and $p$-adic absolute values of its roots. (The $l$-adic absolute values all equal $1$ for $l\neq p$.)

Archimedian estimates have been proved by N. Katz\cite{K,K1,K2}.  We review some of them in Section 2 and state a further conjecture.  Our main focus will be on $p$-adic estimates.  The motivating idea is earlier work of Dwork\cite{Dw1} and Mazur\cite{M1,M2}, who considered zeta functions of smooth projective varieties over ${\mathbb F}_q$.  Specifically, they sought to describe the ($p$-adic) Newton polygons of  polynomial factors of the zeta function.  They accomplished this (Dwork for hypersurfaces, Mazur for complete intersections and more; see also Berthelot and Ogus\cite{BO}) by constructing lower bounds for the Newton polygons from the Hodge numbers of certain related varieties in characteristic zero.  Recently, we gave a new proof of Mazur's theorem for complete intersections using Dwork's original approach(\cite{AS1,AS2}).  The first author's student, John Dollarhide\cite{D}, extended those methods further to give $p$-adic estimates for the roots of $L(t)$ in the case $X={\mathbb P}^n$.  Based on his results, we make some conjectures for more general $X$.  In this setting the (conjectural) lower bounds are constructed from ``Hodge numbers'' of certain logarithmic de Rham complexes on varieties in characteristic zero related to $X$.  Aside from the conjectures, the main new result of this note is that Conjecture~3.1 below is true for $X={\mathbb P}^n$.  The proof is sketched in Section~4.   

Before proceeding further, we give an example to illustrate the form taken by $p$-adic estimates for sums involving multiplicative characters.  We regard the characters $\chi_i$ as taking values in ${\mathbb Q}_p(\zeta_{q-1})$.  The group of characters is cyclic, generated by the Teichm\"uller character $\omega$.  It is defined by the property that for $x\in{\mathbb F}_q^\times$, $\omega(x)$ is the unique $(q-1)$-st root of unity in ${\mathbb Q}_p(\zeta_{q-1})$ which is congruent to $x$ modulo~$p$.  We can write any multiplicative character $\chi$ as a power of $\omega$, say, $\chi = \omega^{-(q-1)e}$ for some $e\in\{0,1/(q-1), \dots,(q-2)/(q-1)\}$.  Consider the Gauss sum
\[ g(\chi,\psi) = \sum_{x\in{\mathbb F}_q}\chi(x)\psi(x), \]
where $\psi:{\mathbb F}_q\to{\mathbb Q}_p(\zeta_p)$ is a nontrivial additive character.  There is an automorphism $\sigma\in{\rm Gal}({\mathbb Q}_p(\zeta_{q-1},\zeta_p)/{\mathbb Q}_p)$ that leaves $\zeta_p$ fixed and sends $\zeta_{q-1}$ to $\zeta_{q-1}^p$.  Since Galois automorphisms leave fixed the $p$-adic valuations of elements of $\bar{\mathbb Q}_p$, the Gauss sums $g(\chi^{p^i},\psi)$ have the same $p$-adic valuation for all $i$.  Write $\chi^{p^i} = \omega^{-(q-1)e^{(i)}}$.  If $q=p^a$ and $\chi$ is nontrivial, then the theorem of Stickelberger says that
\[ {\rm ord}_q\;g(\chi,\psi) = \frac{1}{a}\sum_{i=0}^{a-1} e^{(i)} \]
(${\rm ord}_q$ is the $p$-adic valuation normalized by the condition ${\rm ord}_q\,q = 1$).
Our point is that this type of $p$-adic estimate is typical for sums involving multiplicative characters: the estimate is expressed as the average of certain data associated to each of the conjugates of the sum under the map $\zeta_{q-1}\mapsto\zeta_{q-1}^p$.  

\section{Archimedian estimates}

We recall a result of N. Katz\cite{K} (and we refer to that article for more on the history of archimedian estimates for multiplicative character sums).
Let $D_i\subseteq{\mathbb P}^n$ be the hypersurface $f_i = 0$ and let $D=\bigcup_{i=1}^r D_i$.  We make the following hypotheses:
\begin{equation}
\text{For each $I\subseteq\{1,\dots,r\}$, $X\cap\bigcap_{i\in I} D_i$ is smooth of codimension $|I|$ in $X$.}
\end{equation}
\begin{equation}
\text{Let $b_2,\dots,b_r$ be positive integers prime to $p$ such that $d_1 = \sum_{j=2}^r b_jd_j$.}
\end{equation}
\begin{equation}
\text{The character $\chi^b$ is nontrivial, where $b$ is the l.c.m.\ of $b_2,\dots,b_r$.}
\end{equation}
\begin{theorem}\cite[Theorem~5]{K}
Let $L(t)$ be the $L$-function associated to the sum
\[ \sum_{x\in X({\mathbb F}_q)} \chi(f_1(x))\prod_{j=2}^r \chi^{-b_j}(f_j(x)) = \sum_{x\in X({\mathbb F}_q)} \chi\bigg(\frac{f_1(x)}{f_2(x)^{b_2}\cdots f_r(x)^{b_r}}\bigg). \]
Under the hypotheses $(2.1)$, $(2.2)$, and $(2.3)$, $L(t)^{(-1)^{\dim X+1}}$ is a polynomial of degree equal to the absolute value of the Euler characteristic of $X\setminus D$ and all its reciprocal roots have absolute value $\sqrt{q}^{\,\dim X}$.
\end{theorem}

We conjecture that conditions $(2.2)$ and $(2.3)$ can be weakened.  Consider the exponential sums (1.2).  We make the following hypothesis.
\begin{equation}
\text{The characters $\chi_i$ are all nontrivial.}
\end{equation}
\begin{conjecture}
Let $L(t)$ be the $L$-function associated to the sums $(1.2)$.  \\
{\bf (a)} Under the hypotheses $(1.1)$, $(2.1)$, and $(2.5)$, $L(t)^{(-1)^{\dim X+1}}$ is a polynomial of degree equal to the absolute value of the Euler characteristic of $X\setminus D$ and all its reciprocal roots have absolute value $\sqrt{q}^{\,\dim X}$. \\
{\bf (b)} Under the hypotheses $(1.1)$ and $(2.1)$, if at least one of the characters $\chi_i$ is nontrivial, $L(t)^{(-1)^{\dim X+1}}$ is a polynomial of degree equal to the absolute value of the Euler characteristic of $X\setminus D$ and the absolute values of its reciprocal roots lie in the set $\{\sqrt{q}^{\,i}\mid i=0,1,\dots,\dim X\}$.
\end{conjecture}

Part (a) of the conjecture implies part (b), indeed, one can compute the number of reciprocal roots of $L(t)^{(-1)^{\dim X+1}}$ having a given archimedian size.  For example, suppose $\chi_1$ is the trivial character.  Then
\begin{multline}
\sum_{x\in X({\mathbb F}_q)}\prod_{j=1}^{r}\chi_i(f_i(x)) = \sum_{x\in (X\setminus D_1)({\mathbb F}_q)}\prod_{j=2}^{r}\chi_i(f_i(x)) \\
= \sum_{x\in X({\mathbb F}_q)}\prod_{j=2}^{r}\chi_i(f_i(x)) - \sum_{x\in (X\cap D_1)({\mathbb F}_q)}\prod_{j=2}^{r}\chi_i(f_i(x)).
\end{multline}
If we denote by $L_1(t)$ and $L_2(t)$ the $L$-functions associated to the character sums on the right-hand side, then this equation implies
\[ L(t)^{(-1)^{\dim X+1}} = L_1(t)^{(-1)^{\dim X+1}}L_2(t)^{(-1)^{\dim X\cap D_1 + 1}}. \]
If $\chi_2,\dots,\chi_r$ are all nontrivial, then part (a) of the conjecture applies to both $L$-functions on the right-hand side and one gets a complete archimedian description of $L(t)$.  If some of $\chi_2,\dots,\chi_r$ are trivial, one repeats the process as many times as necessary to express $L(t)$ as a product of $L$-functions that are all described by part~(a) of the conjecture.

\section{$p$-adic estimates}

We consider the case where all $\chi_i$ are nontrivial.  By the discussion of the previous section, the general case can be reduced to this one.  Regard each character as taking values in ${\mathbb Q}_p(\zeta_{q-1})$ and express them as powers of the Teichm\"uller character:
\[ \chi_i = \omega^{-(q-1)e_i},\quad e_i\in\{1/(q-1),\dots,(q-2)/(q-1)\}. \]
Condition (1.1) is equivalent to requiring that $\sum_{i=1}^r e_i d_i$ be an integer.  Let $e=(e_1,\dots,e_r)$ and put
\[ d_e = \sum_{i=1}^r e_i d_i\in{\mathbb Z}. \]
We clearly have $0<d_e<\sum_{i=1}^r d_i$.

Suppose that we were working over ${\mathbb C}$ rather than over ${\mathbb F}_q$ and let $U={\mathbb P}^n\setminus D$.  There is a connection $\nabla_e$ on ${\mathcal O}_U(-d_e)$ defined in homogeneous coordinates as
\[ \nabla_e:=d+\sum_{i=1}^r e_i\frac{df_i}{f_i}\wedge:{\mathcal O}_U(-d_e)\to \Omega^1_U(-d_e). \]
(One can think of this formally as $d:{\mathcal O}_U(-d_e)f_1^{e_1}\cdots f_r^{e_r}\to\Omega^1_U(-d_e)f_1^{e_1}\cdots f_r^{e_r}$ since ``sections'' of ${\mathcal O}_U(-d_e)f_1^{e_1}\cdots f_r^{e_r}$ are formally of degree $0$.)
The pullback of this connection to $X$ has logarithmic poles along $D\cap X$, which is a normal crossing divisor on $X$ by hypothesis (2.1).  One can thus consider the logarithmic subcomplex of this pullback, together with its Hodge filtration (the decreasing filtration obtained by truncating the logarithmic subcomplex on the left).  What is the resulting Hodge filtration on the hypercohomology of this logarithmic subcomplex?

This question is answered in Esnault-Viehweg\cite{EV}.  
Let $\Omega^k_X\langle D\cap X\rangle$ be the module of differential $k$-forms on $X$ with logarithmic poles along $D\cap X$ and let $i:X\hookrightarrow {\mathbb P}^n$ be the inclusion.   Set ${\mathcal L} = i^*{\mathcal O}_{{\mathbb P}^n}(d_e)$ and set $E=\sum_{i=1}^r (q-1)e_i D_i$, an effective divisor on ${\mathbb P}^n$.    Since ${\mathcal O}_{{\mathbb P}^n}(E)\cong {\mathcal O}_{{\mathbb P}^n}(d_e)^{q-1}$, we have ${\mathcal O}_X(E.X)\cong{\mathcal L}^{q-1}$.  By Esnault-Viehweg\cite[Section~(2.7)]{EV}, there is a connection with logarithmic poles
\[ \nabla_e:{\mathcal L}^{-1}\to \Omega^1_X\langle D\cap X\rangle \otimes {\mathcal L}^{-1}, \]
which is just the logarithmic subcomplex of the pullback to $X$ of the connection $\nabla_e$ on $\Omega^\bullet_U(-d_e)$ described above.

The Hodge filtration on the complex $(\Omega^\bullet_X\langle D\cap X\rangle\otimes{\mathcal L}^{-1},\nabla_e)$ gives rise to a spectral sequence
\[ E_1^{ij} = H^j(X,\Omega^i_X\langle D\cap X\rangle\otimes{\mathcal L}^{-1}) \Rightarrow {\mathbb H}^{i+j}(\Omega^\bullet_X\langle D\cap X\rangle\otimes{\mathcal L}^{-1},\nabla_e). \]
By \cite[Section~(2.4)]{EV}, this spectral sequence degenerates at the $E_1$-term.  Futhermore, since $X\setminus D$ is affine and $e_i\not\in{\mathbb Z}$ for all $i$, we have by \cite[Theorem~(2.8)]{EV}
\[ H^j(X,\Omega^i_X\langle D\cap X\rangle\otimes{\mathcal L}^{-1}) = 0\quad\text{for $i+j\neq\dim X$}. \]

Put $h^i(e) = \dim_{\mathbb C}H^{\dim X-i}(X,\Omega^i_X\langle D\cap X\rangle\otimes{\mathcal L}^{-1})$.  We conjecture that
\[ \sum_{i=0}^{\dim X} h^i(e) = |\chi(X\setminus D)|, \]
where $\chi(X\setminus D)$ denotes the Euler characteristic of $X\setminus D$.  When $X={\mathbb P}^n$, this is a consequence of Dollarhide's formulas.  The results of Silvotti\cite{S} provide some heuristic support for this belief in the general case.
Furthermore, if we put $\bar{e} = (1-e_1,\dots,1-e_r)$, then we expect that
\[ h^i(\bar{e}) = h^{\dim X-i}(e). \]
When this is the case, we define $g_e:[0,|\chi(X\setminus D)|]\to[0,\infty)$ to be the function whose graph is the Newton polygon relative to the valuation ${\rm ord}_q$ of the polynomial
\[ \prod_{i=0}^{\dim X} (1-q^it)^{h^i(e)}. \]

We now return to the characteristic $p$ setting.  When one can associate to the variety $X$ and divisor $D$ in characteristic $p$ an appropriate variety and divisor in characteristic zero, one can take the $h^i(e)$ to be the ``Hodge numbers'' of the characteristic zero situation.  A simple case is where $X$ is a smooth complete intersection over ${\mathbb F}_q$.  We take the corresponding variety in characteristic zero to be any smooth complete intersection having the same multidegree.  (So implicitly we are conjecturing that for smooth complete intersections over ${\mathbb C}$, the $h^i(e)$ depend only on $e$, the multidegree, and the $d_i$.)  As in the Gauss sum example in Section~1, we define $e'\in\{1/(q-1),\dots,(q-2)/(q-1)\}^r$ by the condition that $(q-1)e'\equiv p(q-1)e\pmod{q-1}$.  Put $e^{(0)}=e$ and for $i=1,\dots,a-1$, let $e^{(i)} = (e^{(i-1)})'$.  
Let $g:[0,|\chi(X\setminus D)|]\to[0,\infty)$ be the average 
\[ g = \frac{1}{a}\sum_{i=0}^{a-1} g_{e^{(i)}}. \]
\begin{conjecture}
For $X$ a smooth complete intersection in ${\mathbb P}^n$, the Newton polygon of $L(t)^{(-1)^{\dim X+1}}$ relative to the valuation ${\rm ord}_q$ lies on or above the graph of the function $g$.
\end{conjecture}

\section{Dollarhide's results}

Dollarhide's thesis deals with multiplicative character sums on ${\mathbb P}^n$.  Our study of the methods he used and the estimates he obtained led us to the conjectures of the previous section.  In particular, Conjecture~3.1 for $X={\mathbb P}^n$ can be proved by applying what Dwork referred to as ``Laplace transform'' (see \cite[Chapter~11]{Dw}) to the complex studied by Dollarhide.  Dwork's Laplace transform takes Dollarhide's complex to the associated graded of the filtration by order of the pole on the complex $(\Omega^\bullet_{{\mathbb P}^n}(*D)\otimes{\mathcal L}^{-1},\nabla_e)$ (the notation $\Omega^\bullet_{{\mathbb P}^n}(*D)$ indicates we allow poles of arbitrary order along the divisor~$D$).  By \cite{EV}, the inclusion of the logarithmic subcomplex into $\Omega^\bullet_{{\mathbb P}^n}(*D)\otimes{\mathcal L}^{-1}$ (with the filtration by order of the pole) is a filtered quasi-isomorphism.  Furthermore, we speculate that the method of Dollarhide can be generalized to prove Conjecture~3.1 itself.
We give here a brief outline of Dollarhide's work and some examples of the
information contained in his formulas.  

Dwork's theory often allows one to reduce the problem of finding $p$-adic estimates for exponential sums to that of calculating the cohomology of certain complexes in characteristic $p$.  Under condition~(1.1), we have the relation
\begin{equation}
\sum_{x\in{\mathbb A}^{n+1}({\mathbb F}_q)} \chi_1(f_1(x))\cdots \chi_r(f_r(x)) = (q-1)\sum_{x\in{\mathbb P}^n({\mathbb F}_q)} \chi_1(f_1(x))\cdots \chi_r(f_r(x)).
\end{equation}
We are interested in the sum over ${\mathbb P}^n$, but the sum over ${\mathbb A}^{n+1}$ is closely related.  We expect a single nontrivial cohomology group associated to the sum over ${\mathbb P}^n$ in the ``nice'' situation.  The factor of $q-1$ on the right-hand side of~(4.1) tells us that the sum over ${\mathbb A}^{n+1}$ should have associated to it two isomorphic cohomology groups where the Frobenius action on one is simply $q$ times the Frobenius action on the other.  We focus attention on the sum over ${\mathbb A}^{n+1}$.  For a nontrivial multiplicative character $\chi$ and nontrivial additive character $\psi$ one has
\[ \sum_{x\in{\mathbb A}^{n+1}({\mathbb F}_q)} \chi(f(x)) = g(\chi^{-1},\psi)^{-1} \sum_{(x,y)\in{\mathbb A}^{n+2}({\mathbb F}_q)} \chi^{-1}(y)\psi(yf(x)), \]
where $g(\chi^{-1},\psi)$ denotes a Gauss sum (see Section~1).  More generally, 
\begin{multline}
\sum_{x\in{\mathbb A}^{n+1}({\mathbb F}_q)} \chi_1(f_1(x))\cdots \chi_r(f_r(x)) =\\
 \prod_{j=1}^r g(\chi_i^{-1},\psi)^{-1}\sum_{(x,y)\in{\mathbb A}^{n+1+r}({\mathbb F}_q)} \chi_1^{-1}(y_1)\cdots \chi_r^{-1}(y_r)\psi(y_1f_1(x)+\cdots +y_rf_r(x)). 
\end{multline}
Since the $p$-adic valuation of Gauss sums is well known, one is reduced to studying the ``twisted'' exponential sum
\[ \sum_{(x,y)\in{\mathbb A}^{n+1+r}({\mathbb F}_q)} \chi_1^{-1}(y_1)\cdots \chi_r^{-1}(y_r)\psi(y_1f_1(x)+\cdots+ y_rf_r(x)).  \]

Let $\Omega^k_{{\mathbb F}_q[x,y]}\langle Y\rangle$ denote the $k$-forms over ${\mathbb F}_q[x,y]$ with logarithmic poles along the divisor $Y:y_1\cdots y_r=0$ in ${\mathbb A}^{n+1+r}$.  It is a free ${\mathbb F}_q[x,y]$-module with basis the forms ($l+m=k$, $0\leq i_1<\dots<i_l\leq n$, $1\leq j_1<\dots<j_m\leq r$)
\[ dx_{i_1}\cdots dx_{i_l}\frac{dy_{j_1}}{y_{j_1}}\cdots \frac{dy_{j_m}}{y_{j_m}}. \]
Put $F=\sum_{i=1}^r y_if_i(x)\in{\mathbb F}_q[x,y]$ and consider the complex
$(\Omega^\bullet_{{\mathbb F}_q[x,y]}\langle Y\rangle,dF\wedge)$.  By setting
\begin{align*}
\deg_1 x_i &= \deg_1 dx_i = 1&  \deg_1 y_j &= \deg_1 dy_j = -d_j, \\
\deg_2 x_i & = \deg_2 dx_i = 0& \deg_2 y_j &= \deg_2 dy_j = 1,
\end{align*}
this becomes a bigraded complex with boundary map of degree $(0,1)$.  For any $i\in{\mathbb Z}$, we denote by $\Omega^\bullet_{{\mathbb F}_q[x,y]}\langle Y\rangle^{(i)}$ the subcomplex of elements with $\deg_1 = i$, and we denote by $\Omega^k_{{\mathbb F}_q[x,y]}\langle Y\rangle^{(i_1,i_2)}$ the subspace of elements of bidegree $(i_1,i_2)$ in $\Omega^k_{{\mathbb F}_q[x,y]}\langle Y\rangle$.  We denote by
$H^k(\Omega^\bullet_{{\mathbb F}_q[x,y]}\langle Y\rangle^{(i_1)},dF\wedge)^{(i_2)}$ the subspace of elements of bidegree $(i_1,i_2)$ in the $k$-th cohomology group of the complex $(\Omega^\bullet_{{\mathbb F}_q[x,y]}\langle Y\rangle^{(i_1)},dF\wedge)$.

Let $e=(e_1,\dots,e_r)$ be as in Section~3 and set $\bar{e} = (1-e_1,\dots, 1-e_r)$.  We have $d_{\bar{e}} =\sum_{i=1}^r (1-e_i)d_i$ and
$0<d_{\bar{e}}<\sum_{i=1}^r d_i$.  In \cite{D}, Dollarhide (a)~uses Dwork theory to reduce the computation of the lower bound for the Newton polygon of $L(t)^{(-1)^{n+1}}$ to the computation of the cohomology of the complex $(\Omega^\bullet_{{\mathbb F}_q[x,y]}\langle Y\rangle^{(d_{\bar{e}})},dF\wedge)$, and (b)~calculates the $H^k(\Omega^\bullet_{{\mathbb F}_q[x,y]}\langle Y\rangle^{(d_{\bar{e}})},dF\wedge)$.

\begin{theorem}
One has $H^k(\Omega^\bullet_{{\mathbb F}_q[x,y]}\langle Y\rangle^{(d_{\bar{e}})},dF\wedge)=0$ for $k\neq n+r,n+1+r$.
\end{theorem}

Define $\theta:\Omega^k_{{\mathbb F}_q[x,y]}\langle Y\rangle\to\Omega^{k-1}_{{\mathbb F}_q[x,y]}\langle Y\rangle$ to be the ${\mathbb F}_q[x,y]$-module homomorphism acting on basis elements by
\begin{multline}
\theta(dx_{i_1}\cdots dx_{i_l}\frac{dy_{j_1}}{y_{j_1}}\cdots\frac{dy_{j_m}}{y_{j_m}}) = \sum_{s=1}^l (-1)^{s-1}x_s\,dx_{i_1}\cdots \widehat{dx}_{i_s}\cdots dx_{i_l}\frac{dy_{j_1}}{y_{j_1}}\cdots\frac{dy_{j_m}}{y_{j_m}} \\
+(-1)^l\sum_{t=1}^m (-1)^t d_{j_t}\,dx_{i_1}\cdots dx_{i_l}\frac{dy_{j_1}}{y_{j_1}}\cdots\frac{\widehat{dy}_{j_t}}{y_{j_t}}\cdots \frac{dy_{j_m}}{y_{j_m}}. 
\end{multline}
One checks that the sequence
\[ 0\to\Omega^{n+1+r}_{{\mathbb F}_q[x,y]}\langle Y\rangle^{(j)}\xrightarrow{\theta}\Omega^{n+r}_{{\mathbb F}_q[x,y]}\langle Y\rangle^{(j)}\xrightarrow{\theta}\dots\xrightarrow{\theta}\Omega^0_{{\mathbb F}_q[x,y]}\langle Y\rangle^{(j)}\to 0 \]
is exact for all $j>0$.  This leads to the following result.

\begin{corollary}
The map $\theta$ induces an isomorphism
\[ H^{n+r+1}(\Omega^\bullet_{{\mathbb F}_q[x,y]}\langle Y\rangle^{(d_{\bar{e}})},dF\wedge)\cong H^{n+r}(\Omega^\bullet_{{\mathbb F}_q[x,y]}\langle Y\rangle^{(d_{\bar{e}})},dF\wedge). \]
\end{corollary}

In terms of the bigrading, we have a more precise result.
\begin{proposition}
One has
\[ H^{n+r+1}(\Omega^\bullet_{{\mathbb F}_q[x,y]}\langle Y\rangle^{(d_{\bar{e}})},dF\wedge)^{(j)}=0 \quad \text{for $j<0$ or $j>n$} \]
and the map $\theta$ induces isomorphisms
\[ H^{n+r+1}(\Omega^\bullet_{{\mathbb F}_q[x,y]}\langle Y\rangle^{(d_{\bar{e}})},dF\wedge))^{(j)}\cong H^{n+r}(\Omega^\bullet_{{\mathbb F}_q[x,y]}\langle Y\rangle^{(d_{\bar{e}})},dF\wedge)^{(j)} \quad\text{for all $j$.} \]
\end{proposition}

With the computation of cohomology complete, one can apply Dwork theory to find the lower bound for the Newton polygon of $L(t)^{(-1)^{n+1}}$.  For $j=0,\dots,n$, put
\[ k^j(d_{\bar{e}}) = \dim_{{\mathbb F}_q} H^{n+r+1}(\Omega^\bullet_{{\mathbb F}_q[x,y]}\langle Y\rangle^{(d_{\bar{e}})},dF\wedge)^{(j)} \]
and let $\tilde{g}_{d_{\bar{e}}}:[0,|\chi(U)|]\to[0,\infty)$ be the function whose graph is the Newton polygon of $\prod_{j=0}^n (1-q^jt)^{k^j(d_{\bar{e}})}$ relative to the valuation ${\rm ord}_q$.  Put $\tilde{g} = a^{-1}\sum_{i=0}^{a-1} \tilde{g}_{d_{\bar{e}^{(i)}}}$.
The main result of \cite{D} is the following.

\begin{theorem}
If the $f_i=0$ are smooth hypersurfaces in ${\mathbb P}^n$ meeting transversally and all the $\chi_i$ are nontrivial, then $L(t)^{(-1)^{n+1}}$ is a polynomial of degree $|\chi(U)|$ whose Newton polygon relative to the valuation ${\rm ord}_q$ lies on or above the graph of the function $\tilde{g}$.
\end{theorem}

Equation (5.3) below implies that $k^j(d_{\bar{e}}) = k^{n-j}(d_e)$.  Furthermore, it can be shown that the associated graded of the Hodge filtration discussed in Section~3 corresponds to the reverse of the grading by $\deg_2$, i.e., $k^{n-j}(d_e) = h^j(d_e)$.  Thus the function $\tilde{g}$ defined above equals the function $g$ defined in Section~3, so Theorem~4.7 proves Conjecture~3.1 in the case where $X={\mathbb P}^n$.

\section{Explicit formulas}

In this section we give Dollarhide's explicit formulas for the $k^j(d_{\bar{e}})$.  The basic idea is to write down the Hilbert-Poincar\'e series for the $\Omega^k_{{\mathbb F}_q[x,y]}\langle Y\rangle^{(d_{\bar{e}},\bullet)}$ (easily expressible in terms of binomial coefficients) and use the fact that the complex is almost exact to find the Hilbert-Poincar\'e series for 
$H^{n+r+1}(\Omega^\bullet_{{\mathbb F}_q[x,y]}\langle Y\rangle^{(d_{\bar{e}})})^{(\bullet)}$ (which is a polynomial of degree $n$).  We summarize the result of the calculation.

For nonnegative integers $b_1,\dots,b_r$ we write $B=b_1+\cdots+b_r$.  For $B\leq n$, put
\[ a^{(l)}_{b_1\dots b_r} = \frac{(-1)^{n-B}B!}{n!b_1!\cdots b_r!} S_{n-B}(1-l,2-l,\dots,n-l)d_1^{b_1}\cdots d_r^{b_r}, \]
where $S_i$ denotes the $i$-th elementary symmetric function in $n$ variables.  Define a polynomial of degree $n+1$
\[ A_{b_1\dots b_r}(t) = \sum_{l=0}^{n+1}(-1)^l \binom{n+1}{l}a^{(l)}_{b_1\dots b_r}t^l. \]
For any complex number (or indeterminate) $\alpha$ and integer $b\geq 0$, define a polynomial $Q^{(\alpha)}_b(t)\in{\mathbb Q}(\alpha)[t]$ of degree $b$ by the formula
\[ t^{-\alpha}\bigg(t\frac{d}{dt}\bigg)^b\frac{t^\alpha}{1-t} = \frac{Q^{(\alpha)}_b(t)}{(1-t)^{b+1}}. \]
One can show that the polynomial $A_{b_1\dots b_r}(t)$ is divisible by $(1-t)^{B+1}$, so for each vector $\alpha = (\alpha_1,\dots,\alpha_r)$ we can define a polynomial $H_\alpha(t)$ of degree $n$ by the formula
\begin{equation}
H_\alpha(t) = \sum_{b_1+\cdots+ b_r\leq n} A_{b_1\dots b_r}(t)\frac{Q^{(\alpha_1)}_{b_1}(t)\cdots Q^{(\alpha_r)}_{b_r}(t)}{(1-t)^{B+1}}. 
\end{equation}
One can show that $H_\alpha(t)$ depends only on $d_\alpha$, i.e., if $\sum_{i=1}^r \alpha_i'd_i = \sum_{i=1}^r \alpha_id_i$, then $H_{\alpha'}(t) = H_\alpha(t)$.

\begin{proposition}
One has $H_e(t) = \sum_{j=0}^n k^j(d_e) t^j$, the Hilbert-Poincar\'e series of $H^{n+r+1}(\Omega^\bullet_{{\mathbb F}_q[x,y]}\langle Y\rangle^{(d_{e})},dF\wedge)^{(\bullet)}$ for the grading by $\deg_2$.
\end{proposition}

It is clear from the definition that $a^{(n+1-l)}_{b_1\dots b_r} = (-1)^{n-B}a^{(l)}_{b_1\dots b_r}$.  One can show by induction on $b$ that for all $\alpha$ and all $b\geq 0$ one has
\[ t^bQ^{(\alpha)}_b(1/t) = Q_b^{(1-\alpha)}(t). \]
These two relations imply that 
\begin{equation}
H_{\bar{\alpha}}(t) = t^nH_\alpha(1/t) \quad\text{(where $\bar{\alpha} = (1-\alpha_1,\dots,1-\alpha_r)$)}, 
\end{equation}
which implies that $k^j(d_{\bar{e}}) = k^{n-j}(d_e)$.

For small values of $n$, one can use Equation~(5.1) to find explicit formulas for the $k^j(d_e)$.  When $n=1$, one gets 
\[ k^0(d_e) = d_e-1,\quad k^1(d_e) = r-d_e - 1. \]
In the case $n=2$, where the equations $f_i=0$ define smooth curves in ${\mathbb P}^2$, one calculates
\begin{align*}
k^0(d_e) &= \frac{(d_e-1)(d_e-2)}{2} \\
k^1(d_e) &= 1-d_e^2+d_e\sum_{i=1}^r d_i + \frac{1}{2}\sum_{i=1}^r d_i(d_i-3). \\
k^2(d_e) &= \frac{(d_e+1)(d_e+2)}{2} -d_e\sum_{i=1}^r d_i  + \frac{1}{2}\sum_{i=1}^r d_i(d_i-3) + \sum_{1\leq i_1<i_2\leq r} d_{i_1}d_{i_2}.
\end{align*}
For example, take $r=2$, $d_1=d_2$ ($=d$, say), and $e_2=1-e_1$. This implies that $d_e = d$ and the formulas simplify to
\[ k^0(d_e) = k^2(d_e) = \frac{(d-1)(d-2)}{2},\quad k^1(d_e) = 2d^2-3d+1. \]


\begin{thebibliography}{99}

\bibitem{AS1} A. Adolphson and S. Sperber.  On the Jacobian ring of a complete intersection.  J. Algebra {\bf 304} (2006), no.\ 2, 1193--1227.

\bibitem{AS2} A. Adolphson and S. Sperber.  On the zeta function of a projective complete intersection.  Illinois J. Math.\ {\bf 52} (2008), no.\ 2, 389--417.

\bibitem{BO} P. Berthelot and A. Ogus.  Notes on crystalline cohomology.  Princeton University Press, Princeton, N.J.; University of Tokyo Press, Tokyo, 1978.

\bibitem{D} J. Dollarhide. On the $L$-functions of multiplicative character sums. Ph.\ D. dissertation, Oklahoma State University, December, 2008.

\bibitem{Dw1} B. Dwork.  On the zeta function of a hypersurface. II. Ann.\ of Math.\ {\bf 80} (1964), 227--299.

\bibitem{Dw} B. Dwork.  Generalized hypergeometric functions.  Oxford Mathematical Monographs.  Oxford Science Publications.  The Clarendon Press, Oxford University Press, New York, 1990.

\bibitem{EV} H. Esnault and E. Viehweg.  Logarithmic de Rham complexes and vanishing theorems.  Invent.\ Math.\ {\bf 86} (1986), no.\ 1, 161--194.

\bibitem{K} N. Katz. Estimates for nonsingular multiplicative character sums. Int.\ Math.\ Res.\ Not.\ {\bf 2002}, no.\ 7, 333--349.

\bibitem{K1} N. Katz.  Estimates for nonsingular mixed character sums.  Int.\ Math.\ Res.\ Not.\ {\bf 2007}, no.\ 19, Art.\ ID rnm069, 19 pp. 

\bibitem{K2} N. Katz.  Estimates for mixed character sums.  Geom.\ Funct.\ Anal.\ {\bf 18} (2008), no.\ 4, 1251--1269

\bibitem{M1} B. Mazur.  Frobenius and the Hodge filtration.  Bull.\ Amer.\ Math.\ Soc.\ {\bf 78} (1972), 653--667.

\bibitem{M2} B. Mazur. Frobenius and the Hodge filtration (estimates). Ann.\ of Math.\ {\bf 98} (1973), 58--95.

\bibitem{S} R. Silvotti.  On a conjecture of Varchenko.  Invent.\ Math.\ {\bf 126} (1996), no.\ 2, 235--248.
\end{thebibliography}
\end{document}